\documentclass{elsart}

\usepackage{amssymb,amsmath}

\newtheorem{theorem}{Theorem}[section]

\newtheorem{proposition}[theorem]{Proposition}
\newtheorem{corollary}[theorem]{Corollary}

\journal{BULLETIN of the Australian Math. Soc.}
\begin{document}

\begin{frontmatter}

\title{Namioka spaces and topological games}

\author{V.V. Mykhaylyuk}

\address{Chernivtsi National University, Department of Mathematical Analysis,
Kotsjubyns'koho 2, Chernivtsi 58012, Ukraine\\mathan@chnu.cv.ua}

\begin{abstract}
We introduce a class of $\beta-v$-unfavorable spaces, which
contains some known classes of $\beta$-unfavorable spaces for
topological games of Choquet type. It is proved that every
$\beta-v$-unfavorable space $X$ is a Namioka space, that is for
any compact space $Y$ and any separately continuous function
$f:X\times Y\to \mathbb R$ there exists a dense in $X$
$G_{\delta}$-set $A\subseteq X$ such that $f$ is jointly
continuous at each point of $A\times Y$.
\end{abstract}

\begin{keyword}
Namioka space, separately continuous function, Valdivia compact,
topological game

AMS Subject Classification: Primary  54C05, 54D30, 54E52.

\end{keyword}
\end{frontmatter}

\section{Introduction}
\label{Introduction}

Investigations of the joint continuity point set of separately
continuous functions were started by R.~Baire in \cite{B} and were
continued in papers of many mathematicians. A Namioka's result
\cite{N} on a massivity of the joint continuity point set of
separately continuous function on the product of two topological
spaces, one of which satisfies compactness type conditions, gave a
new impulse to a further investigation of this topic.

A topological space $X$ is called {\it a strongly countably
complete space} if there exists a sequence $({\mathcal
U}_n)_{n=1}^{\infty}$ of open coverings of $X$ such that
$\bigcap\limits_{n=1}^{\infty} F_n\ne \O$ for every centered
sequence $(F_n)^{\infty}_{n=1}$ of closed in $X$ sets $F_n$ with
$F_n\subseteq U_n$ for every $n\in \mathbb N$ and some $U_n\in
{\mathcal U}_n$.

\begin{theorem}{\bf (Namioka)}. {\it Let $X$ be a strongly countably
complete space, $Y$ be a compact space and $f:X\times Y \to
\mathbb R$ be a separately continuous function. Then there exists
a dense in $X$ $G_{\delta}$-set $A\subseteq X$ such that $f$ is
jointly continuous at each point of $A\times Y$.}
\end{theorem}

The following notions were introduced in \cite{S}.

A mapping $f:X\times Y \to \mathbb R$ has {\it the Namioka
property} if there exists a dense in $X$ $G_{\delta}$-set
$A\subseteq X$ such that $A\times Y\subseteq C(f)$, where $C(f)$
means the joint continuity point set of $f$.

A Baire space $X$ is called {\it a Namioka space} if for any
compact space $Y$, every separately continuous function $f:X\times
Y\to \mathbb R$ has the Namioka property.

It was shown in \cite{C} that a topological games technique can be
useful in a study of Namioka spaces.

Let ${\mathcal P}$ be a system of subsets of a topological space
$X$. Define a $G_{\mathcal P}$-game on $X$ in which two players
$\alpha$ and $\beta$ participate. A nonempty open in $X$ set $U_0$
is the first move of $\beta$ and a nonempty open in $X$ set
$V_1\subseteq U_0$ and set $P_1\in {\mathcal P}$ are the first
move of $\alpha$. Further $\beta$ chooses a nonempty open in $X$
set $U_1\subseteq V_1$ and $\alpha$ chooses a nonempty open in $X$
set $V_2\subseteq U_1$ and a set $P_2\in {\mathcal P}$ and so on.
The player $\alpha$ wins if
$(\bigcap\limits_{n=1}^{\infty}V_n)\bigcap
(\overline{\bigcup\limits_{n=1}^{\infty}P_n})\ne\O$. Otherwise
$\beta$ wins.

A topological space $X$ is called {\it $\alpha$-favorable in the
$G_{{\mathcal P}}$-game} if $\alpha$ has a winning strategy in
this game. A topological space $X$ is called {\it
$\beta$-unfavorable in the $G_{{\mathcal P}}$-game} if $\beta$ has
no winning strategy in this game. Clearly, any $\alpha$-favorable
topological space $X$ is a $\beta$-unfavorable space.

In the case of ${\mathcal P}=\{X\}$ the game $G_{{\mathcal P}}$ is
the classical Choquet game and $X$ is $\beta$-unfavorable in this
game if and only if $X$ is a Baire space (see \cite{S}). If
${\mathcal P}$ is the system of all finite (or one-point) subsets
of $X$ then $G_{{\mathcal P}}$-game is called a {\it
$\sigma$-game}.

Note that J.~P.~R.~Christensen in \cite{C} generalizing the
Namioka theorem considered an $s$-game which is a modification of
the $\sigma$-game. So, for the $s$-game $\alpha$ and $\beta$ play
in the same way as in the $\sigma$-game, and $\alpha$ wins if each
subsequence $(x_{n_k})_{k=1}^{\infty}$ of sequence
$(x_n)_{n=1}^{\infty}$ has a cluster point in the set
$\bigcap\limits_{n=1}^{\infty} V_n$, where $P_n=\{x_n\}$. It was
proved in \cite{C} that any $\alpha-s$-favorable space is a
Namioka space.

J.~Saint-Raymond showed in \cite{S} that for usage the topological
games method in these investigations it is enough to require a
weaker condition of $\beta$-unfavorability instead of the
$\alpha$-favorability. He proved that any
$\beta-\sigma$-unfavorable space is a Namioka space and
generalized the Christensen result.

A further development of this technique leads to a consideration
of another topological games which based on wider systems
${\mathcal P}$ of subsets of a topological space $X$.

Let $T$ be a topological space and ${\mathcal K}(T)$ be a
collection of all compact subsets of $T$. Then $T$ is said to be
{\it ${\mathcal K}$-countably-determined} if there exist a subset
$S$ of the topological space ${\mathbb N}^{\mathbb N}$ and a
mapping $\varphi:S\to {\mathcal K}(T)$ such that for every open in
$T$ set $U\subseteq T$ the set $\{s\in S: \varphi(s)\subseteq U\}$
is open in $S$ and $T=\bigcup\limits_{s\in S} \varphi (s)$; and it
is called {\it ${\mathcal K}$-analytical} if there exists such a
mapping $\varphi$ for the set $S={\mathbb N}^{\mathbb N}$.

A set $A$ in a topological space $X$ is called {\it bounded} if
for any continuous function $f:X\to \mathbb R$ the set
$f(A)=\{f(a): a\in A\}$ is bounded.

The following theorem gives further generalizations of
Saint-Raymond result.

\begin{theorem}{\it  Any $\beta$-unfavorable in $G_{\mathcal P}$-game topological
space $X$ is a Namioka space if :

$(i)$  ${\mathcal P}$ is the system of all compact subsets of $X$
({\bf M.~Talagrand} \cite{T2});

$(ii)$ ${\mathcal P}$ is the system of all ${\mathcal
K}$-analytical subsets of $X$ ({\bf G.~Debs} \cite{D});

$(iii)$ ${\mathcal P}$ is the system of all bounded subsets of $X$
({\bf O.~Maslyuchenko} \cite{M});

$(iv)$  ${\mathcal P}$ is a system of all ${\mathcal
K}$-countable-determined subsets of $X$ ({\bf V.~Rybakov}
\cite{R}).}
\end{theorem}

It is easy to see that $(iv)\Rightarrow (ii) \Rightarrow (i)$ and
$(iii) \Rightarrow (i)$.

In this paper, using a technique, which is related to the
dependence of functions on products upon some quantity of
coordinates, we prove a result which generalizes $(iii)$ and
$(iv)$.

\section{Dependence of functions upon some
quantity of coordinates and Namioka property}
\label{2}

\begin{proposition} {\it Let $X$ be a topological space,
$A\subseteq X$ be a dense in $X$ set, $Y\subseteq {\mathbb R}^T$
be a topological space, $f:X\times Y\to \mathbb R$ be a continuous
in the first variable function, $\varepsilon \geq 0$ and
$S\subseteq T$ be such that $|f(a,y')-f(a,y'')|\leq \varepsilon$
for every $a\in A$ and every $y',y''\in Y$ with $y'|_S=y''|_S$.
Then $|f(x,y')-f(x,y'')|\leq \varepsilon$ for every $x\in X$ and
every $y',y''\in Y$ with $y'|_S=y''|_S$.}
\end{proposition}

{\bf Proof.} Suppose that $y',y''\in Y$ with $y'|_S=y''|_S$. Put
$h':X\to \mathbb R$, $h'(x)=f(x,y')$, $h'':X\to \mathbb R$,
$h''(x)=f(x,y'')$. Since $h'$ and $h''$ are continuous, the set
$G=\{x\in X: |h'(x)-h''(x)|>\varepsilon\}$ is open. But $G\cap
A=\O$ and $\overline{A}=X$. Thus, $G=\O$ and $f(x,y')=f(x,y'')$
for each $x\in X$.\hfill$\diamondsuit$

\begin{proposition} {\it Let $Y\subseteq {\mathbb R}^T$ be a
compact space, $(Z, |\cdot - \cdot|_Z)$ be a metric space, $f:Y\to
Z$ be a continuous mapping, $\varepsilon \geq 0$ and $S\subseteq
T$ be such that $|f(y')-f(y'')|_Z\leq\varepsilon$ for every
$y',y''\in Y$ with $y'|_S=y''|_S$. Then for every
$\varepsilon'>\varepsilon$ there exist a finite set $S_0\subseteq
S$ and a real $\delta>0$ such that $|f(y')-f(y'')|_Z\leq
\varepsilon'$ for every $y',y''\in Y$ with $|y'(s)-y''(s)|<\delta$
for each $s\in S_0$.}
\end{proposition}

{\bf Proof.} Fix some $\varepsilon'>\varepsilon$. Suppose that the
proposition is false for this $\varepsilon'$. Put
$A=\{(R,n):R\subseteq
S\,\,\,\mbox{is\,\,\,finite\,\,\,and\,\,\,}n\in\mathbb N\}$.
Consider on $A$ the following order: $(R',n')\leq (R'',n'')$ if
$R'\subseteq R''$ and $n'\leq n''$. By the assumption, for every
$a=(R,n)\in A$ there exist $y'_a, y''_a\in Y$ such that
$|f(y'_a)-f(y''_a)|_Z>\varepsilon'$ and
$|y'_a(s)-y''_a(s)|<\frac{1}{n}$ for each $s\in R$. Since $Y^2$ is
a compact then the net $(y'_a,y''_a)_{a\in A}$ has a subnet
$(z'_b,z''_b)_{b\in B}$ which converges  in $Y^2$ to some point
$(y',y'')$. The continuity of $f$ implies
$|f(y')-f(y'')|_Z\geq\varepsilon'>\varepsilon$. For every
$a_0=(\{s_0\},n_0)\in A$ there exists $b_0\in B$ such that $a\geq
a_0$ for every $b\geq b_0$ where $a\in A$ is such that
$(z'_b,z''_b)=(y'_a,y''_a)$. Therefore
$|y'(s_0)-y''(s_0)|\leq\frac{1}{n_0}$ for every $s_0\in S$ and
$n_0\in \mathbb N$. Thus $y'(s)=y''(s)$ for every $s\in S$ and
$|f(y')-f(y'')|_Z\leq\varepsilon$, but it is
impossible.\hfill$\diamondsuit$

\begin{corollary} {\it Let $X\subseteq {\mathbb R}^T$ be a
compact space, $\varphi:X\to {\mathbb R}^S$ be a continuous
mapping and $S_0\subseteq S$ with $|S_0|\leq \aleph_0$. Then there
exists a set $T_0\subseteq T$ with $|T_0|\leq \aleph_0$ such that
$\varphi(x')|_{S_0}=\varphi(x'')|_{S_0}$ for every $x',x''\in X$
with $x'|_{S_0}=x''|_{T_0}$.}
\end{corollary}

{\bf Proof.} Consider a continuous mapping $f:X\to {\mathbb
R}^{S_0}$, $f(x)=\varphi(x)|_{S_0}$. Note that ${\mathbb R}^{S_0}$
is metrizable. Fix a metric $d$ which generates the topology on
${\mathbb R}^{S_0}$. Using Proposition 2.2 for $\varepsilon=0$ we
obtain that for every $n\in\mathbb N$ there exists a finite set
$T_n\subseteq T$ such that $d(f(x'),f(x''))\leq\frac{1}{n}$ for
any $x',x''\in X$ with $x'|_{T_n}=x''|_{T_n}$. The set
$T_0=\bigcup\limits^{\infty}_{n=1}T_n$ is to be
found.\hfill$\diamondsuit$

Let $X$ be a topological space,$(Y,d)$ be a metric space, $f:X\to
Y$ and $A\subseteq X$ be a nonempty set. The real
$\omega_f(A)=\sup\limits_{x',x''\in A}d(f(x'),f(x''))$ is called
{\it the oscillation of $f$ on $A$}, and the real
$\omega_f(x_0)=\inf\limits_{U\in {\mathcal U}}\omega_f(U)$, where
${\mathcal U}$ is the system of all neighborhoods of $x_0\in X$ in
$X$, is called {\it the oscillation of $f$ at $x_0$}.

The following result illustrates relations between the Namioka
property and the dependence of mappings upon some quantity of
coordinates.

\begin{theorem} {\it Let $X$ be a Baire space, $Y\subseteq
{\mathbb R}^T$ be a compact space and $f:X\times Y\to \mathbb R$
be a separately continuous function. Then the following conditions
are equivalent:

$(i)$ $f$ has the Namioka property;

$(ii)$ for every open in $X$ nonempty set $U$ and every
$\varepsilon >0$ there exist an open in $X$ nonempty set
$U_0\subseteq U$ and a set $S_0\subseteq T$ with
$|S_0|\leq\aleph_0$ such that $|f(x,y')-f(x,y'')|\leq \varepsilon$
for every $x\in U_0$ and every $y',y''\in Y$ with
$y'|_{S_0}=y''|_{S_0}$;

$(iii)$ for every open in $X$ nonempty set $U$ and every
$\varepsilon >0$ there exist an open in $X$ nonempty set
$U_0\subseteq U$, a finite set $S_0\subseteq T$ and $\delta
>0$ such that $|f(x,y')-f(x,y'')|\leq \varepsilon$ for every
$x\in U_0$ and every $y',y''\in Y$ with $|y'(s)-y''(s)|<\delta$,
if $s\in S_0$.}
\end{theorem}

{\bf Proof.} $(i)\Rightarrow (ii)$. Let $f$ has the Namioka
property and $U$ be an open in $X$ nonempty set. Then there exists
an open in $X$ nonempty set $U_0\subseteq U$ such that
$|f(x',y)-f(x'',y)|\leq \frac{\varepsilon}{2}$ for every $x', x''
\in U_0$ and every $y\in Y$. Pick an arbitrary point $x_0\in U_0$.
The continuous function $g:Y\to \mathbb R$, $g(y)=f(x_0,y)$
depends upon countable quantity of coordinates, that is there
exists a set $S_0\subseteq T$ with $|S_0|\leq\aleph_0$ such that
$g(y')=g(y'')$ for every $y', y'' \in Y$ with
$y'|_{S_0}=y''|_{S_0}$. For every $x\in U_0$ we have
$|f(x,y')-f(x,y'')|\leq
|f(x,y')-f(x_0,y')|+|f(x_0,y')-f(x_0,y'')|+|f(x_0,y'')-f(x,y'')|\leq
\frac{\varepsilon}{2}+\frac{\varepsilon}{2}=\varepsilon$.

$(ii)\Rightarrow (iii)$. Fix an open in $X$ nonempty set $U$ and
$\varepsilon >0$. By $(ii)$ there exist an open in $X$ nonempty
set $U'\subseteq U$ and a set $S'\subseteq S$ with
$|S'|\leq\aleph_0$ such that $|f(x,y')-f(x,y'')|\leq
\frac{\varepsilon}{2}$ for every $x\in U'$ and every $y',y''\in Y$
with $y'|_{S'}=y''|_{S'}$. Let $S'=\{s_1, s_2, \dots , s_n, \dots
\}$. For every $n\in \mathbb N$ denote by $F_n$ the set of all
points $x\in X$ such that $|f(x,y')-f(x,y'')|\leq \varepsilon$ for
every $y',y''\in Y$ with $|y'(s_k)-y''(s_k)|<\delta$ by $k=1,
2,\dots, n$. The continuity of $f$ in the variable $x$ implies
that all sets $F_n$ are closed. By Proposition 2.2 the continuity
of $f$ in the variable $y$ implies $U'\subseteq
\bigcup\limits_{n=1}^{\infty}F_n$. Since $X$ is a Baire space,
there exist an integer $n_0\in \mathbb N$ and an open in $X$
nonempty set $U_0\subseteq U'$ such that $U_0\subseteq {\rm
int}(F_{n_0})$. It remains to put $\delta=\frac{1}{n_0}$.

$(iii)\Rightarrow (i)$. Suppose that for an open in $X$ nonempty
set $U$ and any real $\varepsilon>0$ there exist an open in $X$
nonempty set $U_0\subseteq U$ and a finite set $S_0\subseteq T$
and a $\delta>0$ such that $|f(x,y')-f(x,y'')|\leq \varepsilon$
for every $x\in U_0$ and every $y',y''\in Y$ with
$|y'(s)-y''(s)|<\delta$, if $s\in S_0$. Then for every
$\varepsilon
>0$ the open set \mbox{$G_{\varepsilon}=\{x\in X: \left(\forall
y\in Y\right) \left(\omega_f(x,y)<\varepsilon\right)\}$} is dense
in $X$. Put $A=\bigcap\limits_{n=1}^{\infty}G_{\frac{1}{n}}$.
Clearly that $f$ is continuous at each point of $A\times Y$. Thus,
$f$ has the Namioka property.\hfill$\diamondsuit$

\begin{proposition} {\it Let $X\subseteq {\mathbb R}^S$ be a
compact space, $f:X\to \mathbb R$ be a continuous function,
$T\subseteq S$ be a set such that $f(x_1)=f(x_2)$ for every
$x_1,x_2\in X$ with $x_1|_T=x_2|_T$, $\varphi:X\to {\mathbb R}^T$,
$\varphi(x)=x|_T$ and $Y=\varphi(X)$. Then the function $g:Y\to
\mathbb R$, $g(x|_T)=f(x)$ is continuous.}
\end{proposition}

{\bf Proof.} Let $y_0\in Y$, $x_0\in X$ be such that
$\varphi(x_0)=y_0$ and $\varepsilon >0$. The set $F=\{x\in X:
|f(x)-f(x_0)|\geq \varepsilon\}$ is closed in $X$, thus $F$ is
compact. Therefore the set $\varphi(F)$ is compact subset of $Y$,
besides $y_0\not\in \varphi(F)$. Thus the set
$V=Y\setminus\varphi(F)$ is a neighborhood of $y_0$. For every
$y\in V$ we have $|g(y)-g(y_0)|<\varepsilon$. Hence $g$ is
continuous at $y_0$.\hfill$\diamondsuit$

\section{Valdivia compacts and the Namioka property}
\label{3}

In this section we establish a result which we will use in the
proof of a generalization of Theorem 1.2.

Recall that a compact space $Y$ is called {\it a Valdivia compact}
if $Y$ is homeomorphic to a compact $Z\subseteq {\mathbb R}^S$
such that a set $B=\{z\in Z:|{\rm supp}\,z|\leq \aleph _0\}$ is
dense in $Z$, where ${\rm supp}\,f$ means {\it the support}
$\{x\in X: f(x)\ne 0\}$ of a function $f:X\to \mathbb R$.

\begin{theorem} {\it Let $X$ be a Baire space,
$Y\subseteq{\mathbb R}^T$ be a Valdivia compact, $\varepsilon
>0$ and $f:X\times Y\to\mathbb R$ be a continuous in the firsts
variable function such that $\omega_{f^x}(y)<\varepsilon$ for
every $x\in X$ and every $y\in Y$, where $f^x:Y\to\mathbb R$,
$f^x(y)=f(x,y)$. Then there exist an open in $X$ nonempty set
$U_0$ and a set $T_0\subseteq T$ with $|T_0|\leq\aleph_0$ such
that $|f(x,y')-f(x,y'')|\leq 3\varepsilon$ for every $x\in U_0$
and every $y',y''\in Y$ with $y'|_{T_0}=y''|_{T_0}$.}
\end{theorem}

{\bf Proof.} Note that $Y$ is homeomorphic to a compact
$Z\subseteq {\mathbb R}^S$ such that the set $B=\{z\in Z:|{\rm
supp}\,z|\leq \aleph _0\}$ is dense in $Z$. Let $\varphi:Y\to Z$
be a homeomorphism and $g:X\times Z\to \mathbb R$,
$g(x,z)=f(x,\varphi^{-1}(z))$. Clearly, $g$ is continuous in the
first variable and $\omega_{g^{x}}(z)<\varepsilon$ for every $x\in
X$ and every $z\in Z$, where $g^x:Z\to\mathbb R$, $g^x(z)=g(x,z)$.

For each $x\in X$ pick a finite covering ${\mathcal W}_x$ of $Z$
by open in $Z$ basic sets such that $\omega_{g^x}(W)<\varepsilon$
for every $W\in {\mathcal W}_x$. For each $W\in {\mathcal W}_x$
choose a finite set $R(W)\subseteq S$ such that for every
$z',z''\in Z$ the conditions $z'\in W$ and $z''(s)=z'(s)$ for any
$s\in R(W)$ imply that $z''\in W$. For a finite set
$S_x=\bigcup\limits_{W\in {\mathcal W}_x}R(W)$ we have
$|g(x,z')-g(x,z'')|<\varepsilon$ for every $z',z''\in Z$ with
$z'|_{S_x}=z''|_{S_x}$. For each $n \in \mathbb N$ we put
$X_n=\{x\in X:|S_x|\leq n\}$. Since $X$ is a Baire space, there
exist an open in $X$ nonempty set $\tilde{U}$ and an integer
$n_0\in \mathbb N$ such that
$\tilde{U}\subseteq\overline{X}_{n_0}$.

Show that there exist an open in $X$ nonempty set $U_0\subseteq
\tilde{U}$ and a set $S_0\subseteq S$ with $|S_0|\leq\aleph_0$
such that $|g(x,b')-g(x,b'')|\leq\varepsilon$ for every $x\in U_0$
and every $b',b''\in B$ with $b'|_{S_0}=b''|_{S_0}$.

Assume that it is false. Pick a set $S_1\subseteq S$ with
$|S_1|\leq\aleph_0$ and an open in $X$ nonempty set
$U_1=\tilde{U}$. By the assumption, there exist $x_1\in U_1$ and
$b_1, c_1\in B$ such that $|g(x_1,b_1)-g(x_1,c_1)|>\varepsilon$
and $b_1|_{S_1}=c_1|_{S_1}$. Using the continuity of $g$ in the
first variable, we find an open in $X$ nonempty set $U_2\subseteq
U_1$ such that $|g(x,b_1)-g(x,c_1)|>\varepsilon$ for every $x\in
U_2$. Put $S_2=S_1\cup({\rm supp}\,b_1)\cup({\rm supp}\,c_1)$. By
the assumption, there exist $x_2\in U_2$ and $b_2,c_2\in B$ such
that $|g(x_2,b_2)-g(x_2,c_2)|>\varepsilon$ and
$b_2|_{S_2}=c_2|_{S_2}$. Doing like that step by step $n_0$ times,
we obtain a decreasing sequence $(U_n)^{n_0+2}_{1}$ of open in $X$
nonempty sets $U_n$, an increasing sequence $(S_n)^{n_0+2}_{n=1}$
of at most countable sets $S_n\subseteq S$ and sequences
$(b_n)^{n_0+1}_{n=1}$ and $(c_n)^{n_0+1}_{n=1}$ of points $b_n,
c_n\in B$ such that the following conditions hold:

$a)$\,\, $U_{n+1}\subseteq U_n$;

$b)$\,\, $S_{n+1}=S_n\cup ({\rm supp}\,b_n)\cup({\rm supp}\,c_n)$;

$c)$\,\, $b_n|_{S_n}=c_n|_{S_n}$;

$d)$\,\, $|g(x,b_n)-g(x,c_n)|>\varepsilon$ for each $x\in U_{n+1}$

for every $n=1,2, \dots ,n_0+1$.

Since $U_{n_0+2}\subseteq U_1=\tilde{U}\subseteq
\overline{X}_{n_0}$, $U_{n_0+2}\cap X_{n_0}\ne \O$. Pick a point
$x_0\in U_{n_0+2}\cap X_{n_0}$ and fix $n\in\{1,2,\dots, n_0+1\}$.
Then $|g(x_0,b_n)-g(x_0,c_n)|>\varepsilon$ by Condition $d)$. The
definition of $S_{x_0}$ implies $b_n|_{S_{x_0}}\ne
c_n|_{S_{x_0}}$. Besides, $b_n|_{S_n}=c_n|_{S_n}$ by Condition
$c)$ and Condition $b)$ implies that $b_n|_{S\setminus
S_{n+1}}=c_n|_{S\setminus S_{n+1}}\equiv 0$. Therefore
$b_n|_{S\setminus (S_{n+1}\setminus S_n)}=c_n|_{S\setminus
(S_{n+1}\setminus S_n)}$. Thus $S_{x_0}\not\subseteq
S\setminus(S_{n+1}\setminus S_n)$, that is
$S_{x_0}\cap(S_{n+1}\setminus S_n)\ne\O$ for every $n=1,2,\dots,
n_0+1$. But this contradicts $|S_{x_0}|\leq n_0$.

Now we show that $|g(x,z')-g(x,z'')|\leq 3\varepsilon$ for every
$x\in U_0$ ³ $z',z''\in Z$ ç $z'|_{S_0}=z''|_{S_0}$.

Fix $x\in U_0$ and $z',z''\in Z$ with $z'|_{S_0}=z''|_{S_0}$.
Since the countably compact set $B$ is dense in $Z$,
$\omega_{g^{x}}(z')<\varepsilon$,
$\omega_{g^{x}}(z'')<\varepsilon$ and $|S_0|\leq\aleph_0$, there
exist $b',b''\in B$ such that $b'|_{S_0}=z'|_{S_0}$,
$b''|_{S_0}=z''|_{S_0}$, $|g(x,z')-g(x,b')|<\varepsilon$ and
$|g(x,z'')-g(x,b'')|<\varepsilon$. Therefore
$b'|_{S_0}=b''|_{S_0}$ and $|g(x,b')-g(x,b'')|<\varepsilon$. Thus
$$
|g(x,z')-g(x,z'')|\leq |g(x,z')-g(x,b')| + |g(x,b')-g(x,b'')| +
|g(x,b'')-g(x,z'')| <
$$
$$
<\varepsilon + \varepsilon +\varepsilon = 3\varepsilon.
$$

Applying Corollary 2.3 to the compact space $Y$, the mapping
$\varphi$ and the set $S_0$, we find an at most countable set
$T_0\subseteq T$ such that $\varphi(y')|_{S_0} =
\varphi(y'')|_{S_0}$ for every $y',y''\in Y$ with
$y'|_{T_0}=y''|_{T_0}$. Then
$|f(x,y')-f(x,y'')|=|g(x,\varphi(y'))-g(x,\varphi(y''))|\leq
3\varepsilon$  for every $x\in U_0$.\hfill$\diamondsuit$

In particular, Theorem 2.4 and Theorem 3.1 imply the following
result which was obtained in \cite[Corollary 1.2]{Bo}.

\begin{corollary} {\it Any separately continuous function on
the product of a Baire space and a Valdivia compact has the
Namioka property, that is any Valdivia compact is a co-Namioka
space.}
\end{corollary}

\section{Namioka spaces and $\beta-v$-unfavorable spaces}
\label{4}

In this section we prove a generalization of Theorem 1.2.

Let ${\mathcal P}$ be a system of subsets of topological space $X$
with the following conditions

$(v_1)$ ${\mathcal P}$ is closed with respect to finite unions;

$(v_2)$ for every set $E\in {\mathcal P}$ and every compact
$Y\subseteq C_p(X)$ a compact $\varphi(Y)$ is a Valdivia compact,
where $\varphi:C_p(X)\to C_p(E)$, $\varphi(y)=y|_E$.

A topological space $X$ is called {\it $\beta-v$-unfaforable} if
the player $\beta$ has no winning strategy in the $G_{{\mathcal
P}}$-game for some system ${\mathcal P}$ with $(v_1)$ and $(v_2)$.

Recall that a compact space $Y$ is called {\it an Eberlein
compact} if $Y$ is homeomorphic to a compact subset of $C_p(X)$
for some compact $X$. A compact space $Y$ is called {\it a Corson
compact} if $Y$ is homeomorphic to a compact $Z\subseteq {\mathbb
R}^S$ such that $|{\rm supp}\,z|\leq \aleph _0$ for every $z\in
Z$. It is known that any Eberlein compact is a Corson compact and
clearly that any Corson compact is a Valdivia compact.

\begin{proposition} {\it Let $X$ be a topological space and
${\mathcal K}$ be a system of all nonempty sets $E\subseteq X$
such that for every compact $Y\subseteq C_p(X)$ a compact
$\varphi(Y)$ is a Corson compact, where $\varphi:C_p(X)\to
C_p(E)$, $\varphi(y)=y|_E$. Then ${\mathcal K}$ has $(v_1)$ and
$(v_2)$.}
\end{proposition}

{\bf Proof.} Let $E_1, E_2\in {\mathcal K}$, $E=E_1\cup E_2$,
$\varphi_1:C_p(X)\to C_p(E_1)$, $\varphi_1(y)=y|_{E_1}$,
$\varphi_2:C_p(X)\to C_p(E_2)$, $\varphi_2(y)=y|_{E_2}$,
$\varphi:C_p(X)\to C_p(E)$, $\varphi(y)=y|_E$ and $Y\subseteq
C_p(X)$ be a compact. Note that a mapping $\psi:\varphi(Y)\to
\varphi_1(Y)\times \varphi_2(Y)$, $\psi(y)=(y|_{E_1}, y|_{E_2})$,
is a homeomorphic embedding. Therefore the compact $\varphi(Y)$ is
a Corson compact. Thus ${\mathcal K}$ has $(v_1)$.

The property $(v_2)$ of system ${\mathcal K}$ is
obvious.\hfill$\diamondsuit$

\begin{proposition} {\it Any $\beta$-unfavorable in
$G_{\mathcal P}$-game topological space $X$, where ${\mathcal P}$
is the system of all bounded subsets of $X$ or ${\mathcal P}$ is a
system of all ${\mathcal K}$-countably-determined subsets of $X$,
is a $\beta-v$-unfavorable space.}
\end{proposition}

{\bf Proof.} Let $E$ be a bounded set in a topological space $X$
and $Y\subseteq C_p(X)$ be a compact. Consider a continuous
mapping $\psi: X\to C_p(Y)$. Clearly, the set $T=\psi(E)$ is
bounded in $C_p(Y)$. Therefore by \cite[Theorem III.4.1]{A} the
closure $\overline{T}$ of $T$ in $C_p(Y)$ is a compact. Then a
compact $Z= \psi_1(Y)$, where $\psi_1:Y\to C_p(T)$,
$\psi_1(y)(t)=t(y)$ for every $y\in Y$ and $t\in T$, is an
Eberlein compact, because $Z$ is homeomorphic to a compact subset
of $C_p(\overline{T})$. Since compacts $Z$ and $\varphi(Y)$, where
$\varphi:C_p(X)\to C_p(E)$, $\varphi(y)=y|_E$, are homeomorphic,
$\varphi(Y)$ is an Eberlein compact, in particular, $\varphi(Y)$
is a Corson compact.

It follows analogously from \cite[Theorem 3.7]{T1} that for every
${\mathcal K}$-countably-determined set $E\subseteq X$ and every
compact $Y\subseteq C_p(X)$ a compact $\varphi(Y)$, where
$\varphi:C_p(X)\to C_p(E)$, $\varphi(y)=y|_E$, is a Corson
compact.

Thus the systems ${\mathcal P}_1$ of all bounded subsets and
${\mathcal P}_2$ of all ${\mathcal K}$-countable-determined
subsets of the topological space $X$ are contained in the system
${\mathcal K}$, by Proposition 4.1. Therefore any
$\beta$-unfavorable space in the $G_{{\mathcal P}_1}$-game or in
the $G_{{\mathcal P}_2}$-game is a $\beta$-unfavorable in
$G_{\mathcal K}$-game and it is $\beta-v$-unfavorable by
Proposition 4.1.\hfill$\diamondsuit$

\begin{theorem} {\it Any $\beta-v$-unfavorable space is a Namioka space.}
\end{theorem}

{\bf Proof.} Let $X$ be a $\beta-v$-unfavorable space. Then there
exists a system ${\mathcal P}$ of subsets $E$ of topological space
$X$ which satisfies $(v_1)$ and $(v_2)$ and such that $X$ is
$\beta$-unfavorable in the $G_{{\mathcal P}}$-game.

Assume that $X$ is not a Namioka space. Then there exist a compact
space $Y$ and a separately continuous function $f:X\times
Y\to\mathbb R$ which does not have the Namioka property. Consider
a continuous mapping $\varphi:Y\to C_p(X)$,
$\varphi(y)(x)=f(x,y)$. Put $Z=\varphi(Y)$. Clearly, $Z$ is a
compact subspace of ${\mathbb R}^X$ and a separately continuous
mapping $g:X\times Z\to \mathbb R$, $g(x,z)=z(x)$, does not have
the Namioka property. Note that $X$ is a $\beta$-unfavorable space
in the Choquet game, i.e. $X$ is a Baire space. Therefore by
Theorem 2.4 there exist an open in $X$ nonempty set $U_0$ and an
$\varepsilon >0$ such that for every open in $X$ nonempty set
$U\subseteq U_0$ and every at most countable set $A\subseteq X$
there exist $x\in U$ and $z',z''\in Z$ such that $z'|_A=z''|_A$
and $|g(x,z')-g(x,z'')|>\varepsilon$.

Show that for every set $E\in {\mathcal P} $ the set $F(E)=\{x\in
U_0:|g(x,z')-g(x,z'')|\leq\frac{\varepsilon}{8}$ {\it for every}
$z',z''\in Z$ {\it with} $z'|_E=z''|_E\}$ is nowhere dense in
$U_0$.

Suppose that it is false. Since $U_0$ is a Baire space with the
topology induced by $X$ and by Proposition 2.1 all sets $F(E)$ are
closed in $U_0$, there exist a set $E_0\in {\mathcal P}$ and an
open in $X$ nonempty set $V_0\subseteq U_0$ such that
$V_0\subseteq F(E_0)$. Suppose $\psi:Z\to {\mathbb R}^{E_0}$,
$\psi(z)=z|_{E_0}$ and $\tilde{Z}=\psi(Z)$. Since the system
${\mathcal P}$ satisfies $(v_2)$, $\tilde{Z}$ is a Valdivia
compact. For every $\tilde{z}\in \tilde{Z}$ choose a point
$\tau(\tilde{z})\in Z$ such that
$\psi(\tau(\tilde{z}))=\tilde{z}$. Consider a mapping $h:V_0\times
\tilde{Z}\to\mathbb R$, $h(x,\tilde{z})=g(x,\tau(\tilde{z}))$.
Since $g$ is continuous in the first variable, $h$ is continuous
in the first variable.

Fix $x_0\in V_0$ and $\tilde{z}_0\in\tilde{Z}$. The set $B=\{z\in
Z:|g(x_0,z)-g(x_0,\tau(\tilde{z}_0))|\geq\frac{\varepsilon}{4}\}$
is a compact subset of $Z$. Besides, $x_0\in F(E_0)$ implies
$\tilde{z}_0\not\in \psi(B)$. Since $\psi$ is continuous, the set
$\psi(B)$ is a compact subset of $\tilde{Z}$. Thus $\psi(B)$ is a
closed subset of $\tilde{Z}$. Therefore the set
$\tilde{W}=\tilde{Z}\setminus\psi(B)$ is a neighborhood of
$\tilde{z}_0$. Then $\tau(\tilde{z})\not\in B$ for every
$\tilde{z}\in\tilde{W}$, that is
$|g(x_0,\tau(\tilde{z}))-g(x_0,\tau(\tilde{z}_0))|=
|h(x_0,\tilde{z})-h(x_0,\tilde{z}_0)|<\frac{\varepsilon}{4}$ for
every $\tilde{z}\in\tilde{W}$. Hence
$\omega_{h^{x_0}}(z_0)\leq\frac{\varepsilon}{4}$, where
$h^{x_0}:\tilde{Z}\to\mathbb R$,
$h^{x_0}(\tilde{z})=h(x_0,\tilde{z})$.

Thus, $h$ satisfies the conditions of Theorem 3.1. Therefore there
exist an open in $X$ nonempty set $\tilde{U}\subseteq V_0$ and an
at most countable set $A_0\subseteq E_0$ such that
$|h(x,\tilde{z}')-h(x,\tilde{z}'')|\leq\frac{3\varepsilon}{4}$ for
every $x\in \tilde{U}$ and every
$\tilde{z}',\tilde{z}''\in\tilde{Z}$ with
$\tilde{z}'|_{A_0}=\tilde{z}''|_{A_0}$.

Pick arbitrary points $x\in \tilde{U}$ and $z',z''\in Z$ such that
$z'|_{A_0}=z''|_{A_0}$. Put $\tilde{z}'=\psi(z')$ and
$\tilde{z}''=\psi(z'')$. Clearly,
$\tilde{z}'|_{A_0}=\tilde{z}''|_{A_0}$. Therefore
$|h(x,\tilde{z}')-h(x,\tilde{z}'')|\leq\frac{3\varepsilon}{4}$.
Since $z'|_{E_0}=\tau(\tilde{z}')|_{E_0}$,
$z''|_{E_0}=\tau(\tilde{z}'')|_{E_0}$ and $x\in\tilde{U}\subseteq
V_0\subseteq F(E_0)$, $|g(x,z')-g(x,\tau(\tilde{z}'))|=
|g(x,z')-h(x,\tilde{z}')|\leq\frac{\varepsilon}{8}$ and
$|g(x,z'')-g(x,\tau(\tilde{z}''))|=
|g(x,z'')-h(x,\tilde{z}'')|\leq\frac{\varepsilon}{8}$. Then
$$
|g(x,z')-g(x,z'')|\leq |g(x,z')-h(x,\tilde{z}')| +
|h(x,\tilde{z}')-h(x,\tilde{z}'')|
$$
$$
+  |h(x,\tilde{z}'')-g(x,z'')| \leq \frac{\varepsilon}{8} +
\frac{3\varepsilon}{4} + \frac{\varepsilon}{8}= \varepsilon.
$$
But this contradicts the choice of $U_0$.

Thus the set $F(E)$ is nowhere dense in $U_0$ for every
$E\in{\mathcal P}$.

Describe a strategy for the player $\beta$ in the $G_{{\mathcal
P}}$-game. The set $U_0$ is the first move of $\beta$. Let
$(V_1,\tilde{E_1})$ be the first move of $\alpha$, where
$V_1\subseteq U_0$ is an open in $X$ nonempty set and
$\tilde{E}_1\in {\mathcal P}$. Then $U_1=V_1\setminus F(E_1)$ is
the second move of $\beta$, where $E_1=\tilde{E}_1$. If
$V_2\subseteq U_1$ is an open in $X$ nonempty set and
$\tilde{E}_2\in {\mathcal P}$ then $U_2=V_2\setminus F(E_2)$ where
$E_2=E_1\cup\tilde{E}_2$. Continuing the procedure of choice by
the obvious manner, we obtain decreasing sequences
$(U_n)^{\infty}_{n=0}$ and $(V_n)^{\infty}_{n=1}$ of open in $X$
nonempty sets $U_n$ and $V_n$ and an increasing sequence
$(E_n)^{\infty}_{n=1}$ of sets $E_n\in {\mathcal P}$ such that
$V_n\subseteq U_{n-1}$, $U_n=V_n\setminus F(E_n)$ and
$\tilde{E}_n\subseteq E_n$ for every $n\in \mathbb N$, where
$\tilde{E}_n\in {\mathcal P}$ is the corresponding part of the
$n$-th move of $\alpha$.

Put $E=\bigcup\limits_{n=1}^{\infty}E_n$. Clearly,
$\bigcup\limits_{n=1}^{\infty}\tilde{E}_n\subseteq E$. Pick a
point $x_0\in\overline{E}$. Note that $g(x_0,z')=g(x_0,z'')$ for
every $z',z''\in Z$ with $z'|_E=z''|_E$, that is the continuous
function $g^{x_0}:Z\to \mathbb R$, $g^{x_0}(z)=g(x_0,z)$, is
concentrated on $E$. Using Proposition 2.2, we obtain that there
exists a finite set $A\subseteq E$ such that
$|g(x_0,z')-g(x_0,z'')|<\frac{\varepsilon}{8}$ for every
$z',z''\in Z$ with $z'|_A=z''|_A$. Pick $n_0\in \mathbb N$ such
that $A\subseteq E_{n_0}$. Then $x_0\in F(E_{n_0})$ therefore
$x_0\not\in U_{n_0}$. Thus $x_0\not\in
\bigcap\limits_{n=0}^{\infty}U_n$ and $\overline{E}\bigcap
(\bigcap\limits_{n=0}^{\infty}U_n)=\O$. In particular,
$\overline{(\bigcup\limits_{n=1}^{\infty}\tilde{E}_n)}
\bigcap(\bigcap\limits_{n=0}^{\infty}U_n)=\O$. Hence the strategy
described above is a winning strategy for $\beta$ in the
$G_{{\mathcal P}}$-game, but it is impossible.

Thus, our assumption is false and the theorem is
proved.\hfill$\diamondsuit$


\begin{thebibliography}{00}

\bibitem {A} {\sc A.~V.~Arhangel'skii}, {\it Topological function spaces}
(Kluwer Acad. Publ., Dordrecht, 1992).

\bibitem  {B} {\sc R.~Baire}, 'Sur les fonctions de variable reelles', {\it An.
Mat. Pura Appl.} (3) (1899) 1-123.

\bibitem {Bo} {\sc A.~Bouziad}, 'Notes sur la propriete de Namioka', {\it Trans.
Amer. Math. Soc.} (3) 344 (1994) 873-883.

\bibitem {C} {\sc J.~P.~R.~Christesen}, 'Joint continuity of separately
continuous functions', {\it Proc. Amer. Math. Soc.} (3) 82 (1981)
455-461.

\bibitem {D} {\sc G.~Debs}, 'Points de continuite d'une function
separement continue', {\it Proc. Amer. Math. Soc.} (1) 97 (1986)
167-176.

\bibitem {M} {\sc O.~V.~Maslyuchenko}, 'Oscillation of separately continuous
functions and topological games', Dys. ... kand. fiz.-mat. nauk.,
Chernivtsi National University, 2002. (in Ukrainian)

\bibitem {N} {\sc I.~Namioka}, 'Separate contimuity and joint
continuity', {\it Pacif. J. Math.} (2) 51 (1974) 515-531.

\bibitem {R} {\sc V.~I.~Rybakov}, 'Some class of Namioka spaces', {\it Mat.
zametki} (2) 73 (2003) 263-268.(in Russian)

\bibitem {S} {\sc J.~Saint-Raymond}, 'Jeux topologiques et espaces de
Namioka', {\it Proc. Amer. Math. Soc.} (3) 87 (1983) 489-504.

\bibitem  {T1} {\sc M.~Talagrand}, 'Espaces de Banach faiblement ${\mathcal
K}$-analytiques', {\it Ann. of Math.} 110 (1979) 407-438.

\bibitem {T2} {\sc M.~Talagrand}, 'Espaces de Baire et espaces de
Namioka', {\it Ann. of. Math.} (2) 270 (1985) 159-164.


\end{thebibliography}
\end{document}